\documentclass[11pt, twocolumn]{article}
\usepackage{authblk}

\usepackage[utf8]{inputenc}
\usepackage[T1]{fontenc}
\usepackage[utf8]{inputenc}
\usepackage{graphicx}
\usepackage{amsmath}
\usepackage{mathtools}
\usepackage{marginnote}
\usepackage[numbers]{natbib}
\usepackage{comment}

\usepackage[normalem]{ulem}

\usepackage{amsmath,amssymb,amsthm}
\usepackage{cases}
\usepackage{paralist}
\usepackage{verbatim}
\usepackage[right]{eurosym} 
\usepackage{url}
\usepackage{enumerate}

\usepackage{booktabs}

\usepackage{mathrsfs}
\usepackage{dsfont}

\usepackage{graphicx}

\usepackage{mathtools}
\usepackage{marginnote}
\usepackage{placeins}
\usepackage{subcaption}

\usepackage{dblfloatfix}

\usepackage[top=2.5cm,bottom=2.5cm,left=2.5cm,right=2.5cm]{geometry}
\usepackage[colorlinks=true, allcolors=blue]{hyperref}

\newtheorem{Def}{Definition}[section]
\theoremstyle{plain}

\newtheorem{Lemma}[Def]{Lemma}
\newtheorem{theorem}[Def]{Theorem}

\newtheorem{Model}[Def]{Model}
\newtheorem{Prop}[Def]{Proposition}

\newtheorem{Prob}[Def]{Control task}

\newcounter{thmpart}[Def]

\newcommand{\thmpart}{%
  \par\refstepcounter{thmpart}%
  \noindent\textup{(\roman{thmpart})} \,  
}

\title{Instantaneous control for stabilization of interacting particle systems: convergence rates and dominant frequencies}

\author{Daniel Jannik Happ\footnote{\href{mailto:dhapp@uni-wuppertal.de}{dhapp@uni-wuppertal.de}} , Birgit Jacob\footnote{ \href{mailto:bjacob@uni-wuppertal.de}{bjacob@uni-wuppertal.de}} ,
Claudia Totzeck\footnote{\href{mailto:totzeck@uni-wuppertal.de}{totzeck@uni-wuppertal.de} }}  
\affil{Port-Hamiltonian Institute, School of Mathematics and Natural Sciences, \\ University of Wuppertal, Germany}
\date{October 2026}

\newcommand{\R}[0]{\mathbb R}
\newcommand{\N}[0]{\mathbb N}

\newcommand{\st}{\text{subject to} \ }

\newcommand{\dx}[1]{\frac{\mathrm{d}}{\mathrm{d}#1}}
\usepackage{mleftright}
\newcommand{\norm}[1]{\lVert #1 \rVert}
\newcommand{\abs}[1]{| #1 |}

\newcommand{\charX}{x}
\newcommand{\charV}{v}
\newcommand{\cmX}{r}
\newcommand{\cmV}{w}
\newcommand{\Psiop}{\Psi}
\newcommand{\Pot}{\mathbb{V}}
\newcommand{\Xeq}{\charX_{\text{des}} }

\newcommand{\meanXeq}{\overline{\charX}_{\text{des}} }
\newcommand{\range}{\mathrm{ran}}
\newcommand{\B}{B_0}
\newcommand{\cB}{\mathbb{B}}
\newcommand{\HS}{E}

\newcommand{\RN}[1]{%
  \textup{\uppercase\expandafter{\romannumeral#1}}%
}
\newcommand{\rn}[1]{%
  \textup{\lowercase\expandafter{\romannumeral#1}}%
}

\usepackage{lipsum}

\begin{document}

\maketitle
\begin{minipage}{0.9\linewidth}
\medskip

\end{minipage}
\begin{abstract}
  We study a class of Cucker--Smale-type systems with a binary interaction potential subject to an external control input. We apply an instantaneous rolling horizon control (IRHC) feedback to address two different control problems: full set-point stabilization of a spatial configuration with zero velocity and stabilization of the center of mass at a desired position with zero mean velocity using a single external agent.
  The IRHC feedback has been shown to preserve the port-Hamiltonian structure of the corresponding open-loop systems. We analyze how the control parameters and the time discretization affect convergence rates and oscillation frequencies and compare the IRHC feedback with the corresponding IDA-PBC feedback obtained by omitting an additional correction term. For the reduced center of mass dynamics, we derive explicit predictions for continuous- and discrete-time spectral decay rates and oscillation frequencies and verify these predictions numerically. For the full nonlinear particle system, the controller comparison is done by numerical experiments. 
\end{abstract}

\section{Introduction}

Interacting particle systems provide a mathematical formulation to describe how binary interactions shape the collective behavior of a group. We consider second-order mechanical systems with Newtonian equations governing each particle's position and velocity coordinates. In \cite{CuckerSmale, CuckerSmalemathem}, Cucker and Smale introduced a velocity alignment term, where each particle accelerates according to a weighted sum of its velocity differences from the other particles. Here, we consider an externally controlled Cucker--Smale-type system with an additional binary interaction potential depending on relative positions. The resulting system admits a port-Hamiltonian formulation \cite{vdSphSsurvey}, where the kinetic energy and the interaction potential define the Hamiltonian, whereas velocity alignment acts as a resistive element \cite{JaTo2024, ErrJaTo2024}. 
A common objective of controlled Cucker--Smale models is velocity alignment with a prescribed common velocity. Different model predictive control approaches have been proposed for this task: In \cite{KaliseRiccati}, a Riccati feedback from a linearized linear-quadratic optimal control problem (see \cite[Chapter $10$]{Trentelman}) was used, whereas in \cite{AlbiPareschi}, a feedback is derived from the first-order optimality conditions of a discretized single time step optimization problem. Requiring the control to be sparse has been studied in \cite{FornasierSparse}. 
Besides velocity alignment, the position of the particles can be included in the control objective. In \cite{CScollisionless}, the relative particle positions are steered to a prescribed value while avoiding interparticle collisions. Instead of assigning relative positions to all particles, the control objective can also include population level quantities. In \cite{InstContrMF}, a small number of external agents is used to steer the center of mass and spatial variance to prescribed values.

In the present work, we study two control tasks: stabilizing the group's center of mass at a prescribed position with zero mean velocity using a single external agent and stabilizing the full configuration by steering each particle to a prescribed position with zero velocity. For the first task, the center of mass and the mean velocity form a reduced subsystem, which we call the \textit{reduced center of mass system}. For the second task, we consider the full particle system, where all particles' positions and velocities are included in the stabilization task.
\\ 
We address these control tasks using \textit{instantaneous rolling horizon control} (IRHC), which is an inexact variant of discrete-time model predictive control \cite{Hinzeinstcontr, Hinzeflcontr}. After discretizing the equations of motion in time, the cost function is evaluated at the current input and at the state predicted after one time step. Instead of solving the resulting optimization problem exactly, at each time step, IRHC computes a single inexpensive gradient update step of an initial control guess. For Cucker--Smale-type systems, \cite{SetPointstab} derives the IRHC set-point feedback using a semi-explicit Euler discretization, identifies its continuous-time closed loop as a port-Hamiltonian system, and establishes sufficient conditions for continuous- and discrete-time convergence.
However, these convergence results do not quantify the speed of convergence and therefore do not describe how it depends on the time step and the gradient step size. This raises the question of how to choose these parameters to achieve rapid convergence while avoiding oscillations around the target. Due to the port-Hamiltonian structure of the underlying open-loop system, the \textit{interconnection and damping assignment passivity-based control} (IDA-PBC) \cite{vdSpHcontrol, IDAPBC} method provides a natural comparison framework. In fact, the IRHC feedback itself admits an IDA-PBC interpretation, but contains an additional correction term. To isolate the effect of this term, we compare the IRHC feedback with the IDA-PBC feedback obtained by omitting it, and study how this difference and the parameter choices affect convergence rates and oscillation frequencies. For the reduced center of mass system, we explicitly compute the spectra of the discrete- and continuous-time system matrices to characterize how the correction term and controller parameters affect these properties. 
For the full nonlinear particle system, we provide numerical experiments to investigate which findings persist beyond the reduced dynamics. These experiments also examine the convergence behavior of both controllers for parameter choices outside the existing sufficient convergence conditions. \\
This paper is structured as follows. In Section \ref{sec:sysclass}, we introduce the Cucker--Smale interacting particle systems and the two control tasks.  Section \ref{sec:IRHC_stability} introduces the IRHC algorithm, the IDA-PBC reference controller and provides the analytical results for the center of mass system. The numerical experiments are presented in Section \ref{sec:numerics}.  

\section{System class and problem formulation}\label{sec:sysclass}
    Consider $N \in\N, \,  N \geq 2$ interacting particles in $\R^d$. Each particle is described by its position $x_i \in \R^d$ and velocity $v_i \in \R^d$ for $i = 1, \ldots, N$. The $N$-particle state space is given by $\HS \times \HS$, where $\HS \coloneq (\R^d)^N$. We collect all positions and velocities in $\charX = (\charX_1, \ldots, \charX_N) \in \HS$ and $\charV = (\charV_1, \ldots, \charV_N) \in \HS$. We equip the space $\HS$ with the inner product   
    \begin{align} 
          \langle \charV, w \rangle &\coloneq \frac{1}{N} \sum_{j= 1}^N v_j \cdot w_j , \,  v, w \in \HS  \, , \label{eq:HSproduct}
    \end{align}
    where $\cdot$ denotes the Euclidean inner product on $\R^d$. 
    The Cucker--Smale models considered here satisfy Newtonian equations of the form 
    \begin{subequations}\label{eq:excharsys}
       \begin{align}
        \dot{\charX} &= \charV , \label{eq:excharsysX} \\ 
         \dot{\charV} &= - \Psiop(\charX) \charV - \nabla \Pot (\charX) + \cB u, \label{eq:excharsysV} \\
        \charX(0) &= \charX_0, \ \charV(0) = \charV_0 \, ,
        \end{align}
    \end{subequations}
    where $\Pot \in C^1(\HS; \R)$ and $\Psiop \colon  \HS \to \R^{Nd \times Nd} = \mathscr{L}(\HS)$. Additionally, $\Psiop(\charX)$ is self-adjoint and positive semidefinite for all $\charX \in \HS$. We use the convention that $\nabla \Pot$ in \eqref{eq:excharsysV} and norms of elements in $\HS$ are computed with respect to the inner product \eqref{eq:HSproduct}.
    The control space $U$ is a finite-dimensional Hilbert space and $\cB \colon U \to \HS$ is linear. 
    We note that every system of the form \eqref{eq:excharsys} admits a port-Hamiltonian formulation. More precisely, \eqref{eq:excharsys} can be written as 
    \begin{align}
        \begin{pmatrix}
            \dot{\charX} \\ 
            \dot{\charV}
        \end{pmatrix}
        &= \left( 
            \begin{pmatrix}
                0&I\\
                -I &0
            \end{pmatrix} 
            - 
            \begin{pmatrix}
                0 &0\\
                0& \Psi(\charX)
            \end{pmatrix}
        \right) \nabla H(\charX, \charV) \nonumber \\
        &\quad + \begin{pmatrix}
            0 \\ 
            \cB 
        \end{pmatrix}
        u, \text{ where } \label{eq:phopenloop} \\
        H(\charX, \charV) &= \frac{1}{2} \norm{\charV}^2 + \Pot(\charX) . \nonumber
    \end{align}
    Subsequently, we introduce the Cucker--Smale realizations of \eqref{eq:excharsys} we study in the present work.  
    \begin{Model}[Fully actuated Cucker--Smale model] \label{Mod:CSparticle}
        The equations of motion for the Cucker--Smale dynamics read
    \begin{subequations}\label{eq:xvparticles}
        \begin{align}
             \dot{x}_i &= v_i, \\
             \dot{v}_i &= \frac{1}{N} \sum_{j=1}^N \psi(\abs{x_j - x_i})(v_j - v_i) \nonumber \\
            &\quad - \frac{1}{N} \sum_{j=1}^N \nabla \mathcal{V} (x_i - x_j) + (\cB u)_i, \label{eq:xvparticlesv} \\
            x_i(0) &=  x_{0,i}, \ v_i(0) =  v_{0,i}, \ i = 1, \ldots, N \, .
        \end{align}
    \end{subequations}
    The first term on the right hand side of \eqref{eq:xvparticlesv} is the Cucker--Smale velocity alignment term. The function $\psi \in C(\R_+;\R_+)$ describes the alignment strength, whereas $\mathcal{V} \in C^1(\R^d; \R)$ denotes the binary interaction potential. We assume that $\nabla \mathcal{V}(-z) = - \nabla \mathcal{V}(z), \ z \in \R^d$, so that the binary interaction forces satisfy Newton's third law. Let the system \eqref{eq:xvparticles} be fully actuated, that is, the control space is given by $U \coloneq \HS$ and the control operator $\cB \in \R^{Nd \times Nd}$ is surjective. Thus, $\range (\cB) = \HS$ and the control can independently accelerate each particle. 
    To cast the system \eqref{eq:xvparticles} into the form \eqref{eq:excharsys}, for $\charX \in \HS$, we define $\Pot(\charX) \in \R$ and $\Psiop(x) \in \R^{Nd \times Nd}$ by 
    \begin{align}
        \Pot(\charX) &= \frac{1}{2 N^2} \sum_{j,k=1}^N \mathcal{V}(x_j - x_k) , \nonumber \\  
        (\Psiop(\charX) \charV)_i &= \frac{1}{N} \sum_{j=1}^N \psi(\abs{x_j - x_i})(v_i - v_j)   \label{eq:psiop} \, .
    \end{align}
    The operator $\Psi(x)$ is self-adjoint and positive semidefinite with respect to the inner product \eqref{eq:HSproduct}.
    \end{Model}
    \begin{Model}[Cucker--Smale model with an external agent]\label{Mod:CSwithAgent}
        Consider a Cucker--Smale model with $\psi, \ \mathcal{V}$ as in Model \ref{Mod:CSparticle}.  Moreover, we consider an external agent who interacts linearly with each particle. Let the position of the external agent be given by $u(t) \in \R^d \eqcolon U$, where $U$ is equipped with the Euclidean inner product. Assume that the acceleration of the $i$th particle towards the position of the agent is given by $\B(u - x_i)$, where $\B \in \R^{d \times d}$ is symmetric positive definite. Then the equations of motion read
        \begin{subequations}\label{eq:xvparticlesAg}
        \begin{align}
             \dot{x}_i &= v_i, \\
             \dot{v}_i &= \frac{1}{N} \sum_{j=1}^N \psi(\abs{x_j - x_i})(v_j - v_i) \nonumber \\
            &\quad - \frac{1}{N} \sum_{j=1}^N \nabla \mathcal{V} (x_i - x_j) + \B( u - \charX_i  ), \label{eq:xvparticlesAgv} \\
            x_i(0) &=  x_{0,i}, \ v_i(0) =  v_{0,i}, \ i = 1, \ldots, N \, .
        \end{align}
    \end{subequations}
        The system \eqref{eq:xvparticlesAg} can be cast into the form \eqref{eq:excharsys}, where $\Psiop$ is given by \eqref{eq:psiop} and $\cB$ is defined as
        \begin{align}
            (\cB u)_j &= \B u , \ u \in \R^d \, . \label{eq:cBwithagent}
        \end{align}
        Since $\B$ is invertible, the range of $\cB$ is given by 
            $\range (\cB) = \{ (y, \ldots, y) \mid y \in \R^d \}$.
        In particular, $\range (\cB)$ is a $d$-dimensional proper subspace of $\HS$ and the input contributes the same acceleration to every particle. The term $- \B \charX_i $ in \eqref{eq:xvparticlesAgv} induces an additional quadratic confinement term in $\Pot$, namely
        \begin{align}
            \Pot(\charX) &= \frac{1}{2 N^2} \sum_{j,k=1}^N  \mathcal{V}(x_j - x_k)  + \frac{1}{2N} \sum_{j=1}^N x_j \hspace{-1pt} \cdot  \hspace{-1pt} \B x_j . \label{eq:Potwithagent} 
        \end{align}
    \end{Model}
    The different ranges of the control operator $\cB$ motivate us to consider different control objectives. To maintain a desired stationary state $(\Xeq, 0) \in \HS \times \HS$, the control must balance the potential forces at $\Xeq$, that is 
    $\nabla \Pot(\Xeq) = \cB u$ for some $u \in U$. For Model \ref{Mod:CSparticle}, this condition is automatically satisfied for all $\Xeq \in \HS$ due to the surjectivity of $\cB$. For Model \ref{Mod:CSwithAgent}, the input forces $(\cB u)_i$ are identical for all $i = 1, \ldots , N$. We will show in Lemma \ref{lemma:meandecomposition} that input operators of this form control the center of mass without changing the relative dynamics. This motivates us to consider full set-point stabilization for Model \ref{Mod:CSparticle} and center of mass stabilization for Model \ref{Mod:CSwithAgent}. 
    \begin{Prob}[Full set-point stabilization]\label{Prob:full}
        Given an admissible desired configuration $\Xeq \in \HS$, that is $\nabla \Pot (\Xeq) \in \range (\cB)$, find $u \colon [0,\infty) \to U$ such that 
        \begin{align*}
            \charX(t) \to \Xeq , \ v(t) \to 0 \text{ as } t \to \infty \, .        \end{align*}
    \end{Prob}
    \begin{Prop}
        For $\charX \in \HS$ and $y \in \R^d$, denote by $\overline{\charX}$ the mean of $\charX$ and by $y \mathds{1}$, the constant vector with value $y$, that is 
        \begin{align*}
            \bar{\charX} \coloneq \frac{1}{N} \sum_{j= 1}^N \charX_j, \
             y \mathds{1} \coloneq (y, \ldots, y ) \, .
        \end{align*}
        Then the Hilbert space $\HS$ splits orthogonally as 
        \begin{align}
            \HS &= \HS_c \oplus \HS_m , \label{eq:Hssplit} \\
            \HS_c &\coloneq \{ y \mathds{1} \mid y \in \R^d \}, \  
            \HS_m \coloneq \{ \charX \in \HS \mid \overline{\charX} = 0 \} \nonumber ,
        \end{align}
        and the corresponding orthogonal projections onto $\HS_c$ and $\HS_m$ are given by 
        \begin{align*}
            P_c \charX = \overline{\charX} \mathds{1}, \ P_m \charX = \charX - \overline{\charX} \mathds{1} \, .
        \end{align*}
    \end{Prop}
    The decomposition \eqref{eq:Hssplit} follows by writing $x = \overline{x} \mathds{1} + (x - \overline{x} \mathds{1})$. The second term has zero mean and therefore is orthogonal to any constant vector. 
    Using the decomposition \eqref{eq:Hssplit}, we can decouple the external agent dynamics \eqref{eq:xvparticlesAg} as follows:
    \begin{Lemma}\label{lemma:meandecomposition}
        Consider Model \ref{Mod:CSwithAgent} and let $\charX, \charV$ solve \eqref{eq:excharsys}.
        Let
        \begin{align*}
            \cmX &\coloneq P_m \charX =   \charX - \overline{\charX} \mathds{1}, \\   
            \cmV &\coloneq P_m \charV = \charV - \overline{\charV} \mathds{1} \, .
        \end{align*}
        Then $\overline{\charX}, \overline{\charV}$ follow the controlled harmonic oscillator dynamics
        \begin{subequations}\label{eq:openloopCOM}
        \begin{align}
            \dot{ \overline{\charX} } &= \overline{\charV}, \\
            \dot{ \overline{\charV} } &= - \B \overline{\charX} + \B u \, ,
        \end{align}
    \end{subequations}
    whereas $\cmX, \cmV$ follow the autonomous dynamics
    \begin{subequations}\label{eq:uncontcharsys}
       \begin{align}
        \dot{\cmX} &= \cmV , \label{eq:uncontcharsysX} \\ 
         \dot{\cmV} &= - \Psiop(\cmX) \cmV - \nabla \Pot (\cmX) , \label{eq:uncontcharsysV} \\
        \cmX (0) &= \charX_0 - \overline{\charX}_0 \mathds{1} , \ \cmV (0) = \charV_0 - \overline{\charV}_0 \mathds{1} \, .
        \end{align}
    \end{subequations}
    \end{Lemma} 
    For Model \ref{Mod:CSwithAgent}, the decoupling into \eqref{eq:openloopCOM} and \eqref{eq:uncontcharsys} follows from the fact that the alignment and binary interaction forces have zero mean.
    \begin{Prob}[Center of mass stabilization]\label{Prob:mean}
        Given $\meanXeq \in \R^d$, find $u \colon [0,\infty) \to U$ such that 
        \begin{align*}
            \bar{\charX}(t) \to \meanXeq , \ \overline{\charV}(t) \to 0 \text{ as } t \to \infty \, .
        \end{align*}
    \end{Prob}

    We note that $(\meanXeq,0) \in \R^d \times \R^d $ is always an admissible equilibrium for the center of mass dynamics \eqref{eq:openloopCOM} since choosing $u = \meanXeq$ makes the right hand side of \eqref{eq:openloopCOM} vanish at $(\overline{\charX}, \overline{\charV}) = (\meanXeq,0)$.

   \section{IRHC, IDA-PBC comparison and their stability properties}\label{sec:IRHC_stability}
   \subsection{Description of the algorithm}\label{subsec:IRHCdesc}
   Both control tasks \ref{Prob:full} and \ref{Prob:mean} can be addressed using the IRHC approach, which we introduce here. Starting from the current state, we predict the system over a single time step by choosing a numerical integration method.
   For a given control input, we then evaluate a quadratic cost function at the predicted state. Instead of minimizing the cost exactly, we perform a single gradient update of the control. 
   The next state of the system is then computed using the updated control and the procedure is repeated.
    In the following, we apply the IRHC approach to the Cucker--Smale model. For a time-step $\Delta t > 0$, we discretize \eqref{eq:excharsys} using the semi-explicit Euler scheme
   \begin{subequations}\label{eq:semexupd}
   \begin{align}
       \charX^{l+1} &= \charX^l + \Delta t \charV^{l+1}, \label{eq:semexXupd} \\ 
       \charV^{l+1} &= \charV^{l} 
       - \Delta t \big(\Psiop(\charX^l) \charV^{l} + \nabla \Pot (\charX^l) \big) 
       + \Delta t \cB u^l,  \label{eq:semexVupd} \\
        \charX^0 &= \charX_{0}, \quad \charV^0 = \charV_{0}. 
   \end{align}
   \end{subequations}
   Fix $c > 0$. Given an iteration $l \in \N_0$ with current state $(\charX^l , \charV^l)$, consider the minimization problem
   \begin{align}
       \min_{u^l  \in U } J(\charX, \charV, u^l) &\coloneq \Delta t H_{\text{des}}(\charX, \charV ) \label{eq:costJ} \\ 
       \st  \charX &= \charX^{l+1}, \ \charV = \charV^{l+1}  \text{ as in } \eqref{eq:semexXupd}, \eqref{eq:semexVupd} , \nonumber
   \end{align}
   where the desired Hamiltonian is given by 
   \begin{align}
       H_{\text{des}}(x,v) = \frac{1}{2} \norm{\charV }^2 + \frac{c}{2} \norm{\charX - \Xeq}^2 \, . \label{eq:Hes}
   \end{align}
   Since $\charX^{l+1} $ and $\charV^{l+1}$ are uniquely determined by $u^l$, we can define the reduced cost function $f^l$ by $f^l(u^l) \coloneq J( \charX^{l+1}, \charV^{l+1}, u^l)$, where $(\charX^{l+1},\charV^{l+1})$ denote the states obtained from \eqref{eq:semexXupd} and \eqref{eq:semexVupd} for the control $u^l$. Then a single gradient update of an initial control $u_0^l \in U$ with step-size $\rho > 0$ yields the control update
   \begin{align}
       u^l = u_0^l  - \rho \nabla f^l(u_0^l) \, . \label{eq:Gradupdcontrol}
   \end{align}
   For the center of mass stabilization problem \ref{Prob:mean}, we identify $\meanXeq \in \R^d$ with the constant vector $\Xeq = \meanXeq \mathds{1}$. Then due to the orthogonal decomposition $\HS = \HS_c + \HS_m$, we have for the predicted states $\charX^{l+1}, \charV^{l+1} \in \HS$
   \begin{align*}
       \norm{\charX^{l+1} - \Xeq}^2 &= \abs{\overline{\charX}^{l+1} - \meanXeq}^2 + \norm{P_m \charX^{l+1}}^2, \\ 
       \norm{\charV^{l+1} }^2 &= \abs{\overline{\charV}^{l+1} }^2 + \norm{P_m \charV^{l+1} }^2 .
   \end{align*}
   According to Lemma \ref{lemma:meandecomposition}, the centered position and velocity $P_m \charX^{l+1}, P_m \charV^{l+1}$ are independent of the control. Therefore, they do not affect the gradient of the reduced cost with respect to $u$ and the controller obtained from the full dynamics and $H_{\text{des}}$ is the same as the one obtained by applying the IRHC method to the center of mass dynamics \eqref{eq:openloopCOM}. 
   In the present work, we always use a constant initial control $u_0^l = \hat{\alpha}$, where $\hat{\alpha}$ makes the desired configuration an equilibrium, that is $\cB \hat{\alpha} = \nabla \Pot (\Xeq)$. \\
    It was shown in \cite{SetPointstab} that the IRHC algorithm as described above preserves the port-Hamiltonian structure \eqref{eq:phopenloop} of the open-loop system \eqref{eq:excharsys}. In the following, we summarize the main results: 
    \begin{Lemma}{\cite[Lemma $4.2$]{SetPointstab}}\label{lemma:IDA-PBC}
        Let $c_1,c_2 \geq 0$ and consider the system \eqref{eq:excharsys} with the feedback
        \begin{align}
            \alpha (\charX, \charV) &= - c_1 \cB^* \charV + c_2 \cB^* ( \Psiop(\charX) \charV + \nabla \Pot(\charX) ) \nonumber \\
            &\quad - F(\charX), \label{eq:Feedb}
        \end{align}
        where $F \colon \HS \to U$. Suppose that $c_2 \norm{\cB \cB^*} < 1$ and that there exists $\hat{\Pot} \in C^1(\HS ; \R)$ such that $B F(\charX) = T \nabla \hat{\Pot} (\charX), \charX \in \HS$, where $T \coloneq I - c_2 \cB \cB^*$.
        Then the closed-loop system with feedback \eqref{eq:Feedb} has the port-Hamiltonian formulation
        \begin{subequations}\label{eq:pHclosedloop}
            \begin{align}
                 \begin{pmatrix}
                \dot{\charX} \\ 
                \dot{\charV} 
            \end{pmatrix}
            \hspace{-2pt}
            &= \hspace{-1pt} \left( \hspace{-2pt}
            \begin{pmatrix}
                0 & T \\
                -T & 0
            \end{pmatrix}
            \hspace{-1pt}
            - \begin{pmatrix}
                0 & 0\\
                0&  R(\charX)
            \end{pmatrix}
            \hspace{-2pt}
             \right) 
             \hspace{-2pt}
            \nabla H_{\text{cl}}(\charX, \charV) , \\ 
            \charX (0) &= \charX_0 , \ \charV (0) = \charV_0 ,
            \end{align}
        \end{subequations}
        where 
        \begin{align}
            R(\charX) &= T \Psiop(\charX) T + c_1 T^{\frac{1}{2}} \cB \cB^* T^{\frac{1}{2}}, \label{eq:Rcl} \\
            H_{\text{cl}}(\charX, \charV) &= \frac{1}{2} \langle \charV, T^{-1} \charV \rangle + \Pot(\charX) + \hat{\Pot}(\charX) \nonumber .
        \end{align}
    \end{Lemma}
    \begin{theorem}[Port-Hamiltonian formulation {\cite[Corollary~4.6]{SetPointstab}}] \label{Thm:IRHCpHclosedloop}
        We consider the instantaneous control algorithm with cost function \eqref{eq:costJ}, semi-explicit Euler scheme \eqref{eq:semexupd} with step-size $\Delta t > 0$ and gradient step-size $\rho > 0$. Define 
        \begin{align}
            c_1 &\coloneq \rho (\Delta t)^2 \big(1+c(\Delta t)^2 \big) , \nonumber \\
            c_2 &\coloneq \rho (\Delta t)^3 \big(1+c(\Delta t)^2 \big), \
            c_3 \coloneq c \rho (\Delta t)^3 \, . \label{eq:cContIRHC}
        \end{align}
        Suppose that there exists $\hat{\alpha} \in U$ such that $\cB \hat{\alpha} = \nabla \Pot( \Xeq) $ and choose $u_0^l = \hat{\alpha}$. Then the IRHC method computes the control 
        \begin{align}
            u^l &= - c_1 \cB^* \charV^l + c_2 \cB^* ( \Psiop(\charX^l) \charV^l + \nabla \Pot(\charX^l) ) \nonumber \\
            &\quad - c_3 \cB^* (\charX^l - \Xeq) - c_2 \cB^* \nabla \Pot (\Xeq) + \hat{\alpha} . \label{eq:IRHCcontrol}
        \end{align}
        In particular, if $c_2 \norm{\cB^* \cB} < 1$, then the IRHC algorithm computes the semi-explicit Euler iterates with step-size $\Delta t$ of the system \eqref{eq:pHclosedloop} with
        \begin{align}
        \hat{\Pot} (\charX)  &= \frac{c_3}{2} \norm{ \cB^* T^{-\frac{1}{2}}(\charX - \Xeq) }^2 \nonumber \\
            &\quad - \langle \nabla \Pot(\Xeq), \charX - \Xeq \rangle \, . \label{eq:hatPotstatad}
        \end{align}
            
    \end{theorem}

    \subsection{Controller comparison and convergence baseline}\label{subsec:ConrComp}
    In the previous subsection, we have seen that the IRHC method computes the control \eqref{eq:IRHCcontrol}, whereas Lemma \ref{lemma:IDA-PBC} provides a sufficient condition to ensure that a feedback of the form \eqref{eq:Feedb} yields a port-Hamiltonian closed-loop system. In particular, we note that Lemma \ref{lemma:IDA-PBC} does not require $c_2 > 0$. Thus for $c_2 = 0$ and 
    \begin{align*}
        F_0(\charX) \coloneq c_3 \cB^* (\charX - \Xeq) - \hat{\alpha},
    \end{align*}
     we have $T = I$ and 
     \begin{align*}
         &\cB F_0(\charX) = \nabla \Pot_0(\charX), \text{ where }\\
         &\Pot_0(\charX) \hspace{-1.5pt} = \hspace{-1pt} \frac{c_3}{2} \norm{\cB^* (\charX - \Xeq)}^2 \hspace{-2pt} -  \hspace{-1pt} \langle \nabla \Pot (\Xeq) , \charX - \Xeq \rangle . 
     \end{align*}
     Therefore, Lemma \ref{lemma:IDA-PBC} establishes that 
     \begin{align}
         u_{\text{PBC}} = - c_1 \cB^* \charV - c_3 \cB^* (\charX - \Xeq) + \hat{\alpha} \label{eq:IDAPBCcontrol}
     \end{align}
     is an IDA-PBC feedback. In contrast, the IRHC controller \eqref{eq:IRHCcontrol} computes a similar feedback which contains the additional $c_2$ term. To investigate the effect of this term, it is natural to keep the constants $c_1$ and $c_3$ as in \eqref{eq:cContIRHC} while comparing the controls \eqref{eq:IRHCcontrol} and \eqref{eq:IDAPBCcontrol}. 
     In addition, although Theorem \ref{Thm:IRHCpHclosedloop} establishes a port-Hamiltonian closed-loop structure, we have not yet discussed the question of how to choose the parameters $\Delta t, \rho$ and $c$. In order to simplify the notation, when $c_1,c_2,c_3$ are given by \eqref{eq:cContIRHC}, we define $\eta \coloneq c_2 \norm{\cB \cB^*}$. Then we have 
     \begin{align}
         c_1 &=  \frac{1}{\Delta t \norm{\cB \cB^*}} \eta, \
         c_2 = \frac{1}{\norm{\cB \cB^*}} \eta, \nonumber \\
         c_3 &= \frac{c}{(1 + c (\Delta t)^2) \norm{\cB \cB^*}  } \eta \label{eq:cetaPar}
     \end{align}
     and the condition $c_2 \norm{\cB \cB^*} < 1$ becomes $\eta < 1$.
     In the following, we will distinguish between different feedback laws and, given a fixed feedback law, different time evolutions. For fixed parameters $c_1,c_2,c_3$, we refer to the port-Hamiltonian ordinary differential equation \eqref{eq:pHclosedloop} with the corresponding fixed coefficients $c_1,c_2,c_3$ as the \textit{continuous-time closed loop}. Similarly, we call its semi-explicit Euler discretization with step-size $\Delta t > 0$ the \textit{discrete-time closed loop}. 
     Fix $ 0 \leq \eta < 1, \Delta t , c > 0$ and consider the feedback controller \eqref{eq:Feedb} with 
     \begin{align*}
         F(\charX) = c_3 \cB^* (\charX - \Xeq) + c_2 \cB^* \nabla \Pot (\Xeq) - \hat{\alpha} .
     \end{align*}
     If the parameters $c_1  ,c_2  ,c_3$ are given by \eqref{eq:cetaPar}, we call the corresponding feedback law the \textit{IRHC controller}. If $c_1, c_3$ are given by \eqref{eq:cetaPar} and $c_2 = 0$, we call the corresponding feedback law the \textit{IDA-PBC controller}. 
     In particular, the iterates of the IRHC algorithm coincide with the iterates of the discrete-time closed loop system using the IRHC controller. Thus, in continuous-time, $\Delta t$ enters only through the controller parameters \eqref{eq:cetaPar}, whereas in discrete-time, the same $\Delta t$ is additionally used as the time step of the semi-explicit Euler scheme.
     A first baseline for the parameter choice is provided by the following theorem:
     \begin{theorem}[Analytical convergence baseline for the IRHC controller {\cite[Proposition~3.3, Theorem~5.12 and Corollary~5.15]{SetPointstab}} ]\label{thm:convbaseline}
    Consider the full set-point stabilization problem \ref{Prob:full} for Model \ref{Mod:CSparticle}, assume that $\cB$ is surjective,  and that $\nabla \mathcal{V} \colon \R^d \to \R^d $ is Lipschitz continuous with Lipschitz constant $L$. Let $\psi \in C(\R_+, \R_+) $ be bounded and assume that $\psi(x) \geq \underline{\psi}$ for some $\underline{\psi} > 0$. Then there exists $\sigma > 0$ such that 
    \begin{align*}
        \norm{\cB^* \charX} \geq \sigma \norm{\charX} \text{ for all } \charX \in \HS .
    \end{align*}
    Define $K \coloneq 4 L / {\sigma^2}$. We parametrize $c_1,c_2,c_3$ by $\eta, c$ and $\Delta t$ as in \eqref{eq:cetaPar}, where $\eta \coloneq c_2 \norm{\cB \cB^*}$.
    \thmpart\label{thm:baselineconttimes} Continuous time convergence: If 
        \begin{align}
            K \Big( (\Delta t)^2 + \frac{1}{c} \Big) \norm{\cB \cB^*} &< \eta < 1, \label{eq:cteta} \\
            0 < \Delta t < \frac{1}{\sqrt{K \norm{\cB \cB^*} } }, 
            c &> \frac{K \norm{\cB \cB^*}}{ 1 - K \norm{\cB \cB^*} (\Delta t)^2}, \label{eq:ctdtc}
        \end{align}
        then for any $\charX_0, \charV_0 \in \HS$, the continuous-time closed loop solution satisfies 
        $\charX(t) \to \Xeq$ and $\charV (t) \to 0$ as $t \to \infty$. 
    \thmpart\label{thm:baselinediscrtimes} Discrete time convergence: If $ 0 < \eta \leq 1/{66}$ and 
    \begin{align}
        c > \frac{K \norm{\cB \cB^*}}{\eta},
    \end{align}
     then there exists $h > 0$ such that for any initial condition $\charX_0, \charV_0 \in \HS$ and any step-size $0 < \Delta t \leq h$, the iterates of the IRHC algorithm satisfy 
     $\charX^l \to \Xeq$ and $\charV^l \to 0$ as $l \to \infty$. 
    \end{theorem}
   The proof of Theorem \ref{thm:convbaseline} is technical. The main tool for the discrete-time convergence is an energy balance for a modified Hamiltonian. Here, choosing $\eta \leq 1/{66}$ in combination with the conditions stated in Theorem \ref{thm:baselinediscrtimes} is sufficient for this energy balance to be satisfied. We will explore the effect of the parameters in more detail in Sections \ref{subsec:comstability} and \ref{sec:numerics}.
   
   \subsection{Parameter implications of the center of mass dynamics}\label{subsec:comstability}

   Since the equations of motion for the center of mass \eqref{eq:openloopCOM} are linear, finite dimensional and independent of the number of particles, stability, decay rates and oscillation frequencies of the closed-loop system can be computed explicitly as the following lemma shows:
   \begin{Lemma}\label{lemma:COMspectrum}
       For the external agent Model \ref{Mod:CSwithAgent}, consider the closed-loop system \eqref{eq:pHclosedloop}, where $\hat{\Pot}$ is given by \eqref{eq:hatPotstatad}. Then the shifted center of mass coordinates $\overline{\charX} - \meanXeq , \overline{\charV}$ satisfy 
       \begin{align}
            &\dx{t} 
            \begin{pmatrix}
                \overline{\charX} - \meanXeq \\ 
                \overline{\charV}
            \end{pmatrix}
            = A_{ \text{ct} }
            \begin{pmatrix}
                \overline{\charX} - \meanXeq \\ 
                \overline{\charV}
            \end{pmatrix}, \label{eq:CTreduced} \\ 
            &\text{where } A_{ \text{ct} } \coloneq \begin{pmatrix}
                0 & I \\ 
                -( \B + c_3 \B^2 - c_2 \B^3) & - c_1 \B^2 
            \end{pmatrix}
            \nonumber .
       \end{align}
       Moreover, the corresponding semi-explicit Euler discretization reads 
       \begin{align}
           &\begin{pmatrix}
                \overline{\charX}^{l+1} - \meanXeq \\ 
                \overline{\charV}^{l+1}
            \end{pmatrix}
            = A_{ \text{dt} }
            \begin{pmatrix}
                \overline{\charX}^l - \meanXeq \\ 
                \overline{\charV}^l
            \end{pmatrix}, \label{eq:DTreduced} \\ 
             &A_{ \text{dt} } \coloneq I + \Delta t \begin{pmatrix}
                 I &\Delta t I \\
                 0& I 
             \end{pmatrix}
             A_{ct} . \nonumber
       \end{align}
       For $\lambda \in \sigma (\B)$, we define the $2 \times 2$ matrices
       \begin{align}
            A_{\text{ct}}(\lambda) &\coloneq 
            \begin{pmatrix}
                0 & 1 \\ 
                -( \lambda + c_3 \lambda^2 - c_2 \lambda^3) & - c_1 \lambda^2 
            \end{pmatrix} , \label{eq:reprmatrix} \\ 
            A_{\text{dt}}(\lambda) &\coloneq I 
            + \Delta t 
            \begin{pmatrix}
                1 &\Delta t \\ 
                0 & 1
            \end{pmatrix}
            A_{\text{ct}}(\lambda) \nonumber . 
        \end{align}
        Then their eigenvalues are given by
        \begin{align}
           &\mu_{\text{ct}, \pm }(\lambda) \nonumber \\
           &= \frac{ \mathrm{tr} (A_{\text{ct} } (\lambda ) ) }{2} \pm 
           \sqrt{ \Big( \frac{ \mathrm{tr} (A_{\text{ct} }(\lambda  ) ) }{2}\Big)^2 - \det (A_{\text{ct} } (\lambda ) )  } \nonumber \\ 
             &=
            - \frac{c_1 \lambda^2}{2} 
            \pm \sqrt{ \Big( \frac{c_1 \lambda^2}{ 2 } \Big)^2 - ( \lambda + c_3 \lambda^2 - c_2 \lambda^3 ) }, \label{eq:COMlambdaCT} \\
           &\mu_{\text{dt}, \pm }(\lambda) \nonumber \\
           &= \frac{ \mathrm{tr} (A_{\text{dt} }(\lambda) ) }{2} \pm 
           \sqrt{ \Big( \frac{ \mathrm{tr} (A_{\text{dt} }(\lambda) ) }{2}\Big)^2 - \det (A_{\text{dt}} (\lambda) )  } ,   \\
           &\mathrm{tr} (A_{ \text{dt} }  (\lambda)) = 2 - c_1 \Delta t \lambda^2  \nonumber \\
           &\hspace{2cm} - (\Delta t)^2 (\lambda + c_3 \lambda^2 - c_2 \lambda^3) , \label{eq:COMlambdaDT} \\ 
           &\det (A_{\text{dt} } (\lambda) ) =  1 - c_1 \Delta t \lambda^2 . \nonumber
       \end{align}
       Then the spectra of $A_{ \text{ct}}$ and $A_{\text{dt} }$ are given by 
       \begin{align}
           \sigma ( A_{ \text{ct} } ) &= \{ \mu_{\text{ct}, \pm} (\lambda) \mid \lambda \in \sigma(\B) \} , \nonumber \\
           \sigma (A_{ \text{dt} } ) &= 
           \{ \mu_{\text{dt}, \pm} (\lambda) \mid \lambda \in \sigma(\B) \} . \label{eq:sigmadt}
       \end{align}
   \end{Lemma}
   \begin{proof}
       The continuous- and discrete-time equations of motion for the reduced center of mass dynamics \eqref{eq:CTreduced} and \eqref{eq:DTreduced} follow by taking the mean in \eqref{eq:pHclosedloop} and its semi-explicit Euler discretization and by substituting the expressions \eqref{eq:cBwithagent} and \eqref{eq:Potwithagent} for $\cB$ and $\Pot$. \\
       If $q \in \R^d \setminus \{0 \}$ is an eigenvector of $\B$ with eigenvalue $\lambda$, then the two-dimensional subspace 
       \begin{align*}
           V_q \coloneq \mathrm{span} \Big\{ \begin{pmatrix}
               q \\ 0
           \end{pmatrix} 
           ,
           \begin{pmatrix}
               0 \\ q
           \end{pmatrix} 
           \Big\} \subseteq \R^d \times \R^d
       \end{align*}
       is invariant under $A_{\text{ct}}$ and its restriction to $V_q$ with respect to the basis $(q,0),(0,q)$ has the representation $A_{\text{ct}}(\lambda)$.
       A short computation shows that the spectrum of the $2 \times 2$ matrix \eqref{eq:reprmatrix} is given by $\mu_{\text{ct}, \pm}(\lambda)$. Since $\B$ is symmetric, its eigenspaces span all of $\R^d$ and therefore taking the union over all $\lambda \in \sigma (\B)$ yields the spectrum of $A_{\text{ct}}$. The subspace $V_q$ is also invariant under $A_{\text{dt}}$ and thus a similar argument yields \eqref{eq:sigmadt}. 
   \end{proof}

\begin{figure}[ht]
        \centering
            \includegraphics[width=\linewidth]{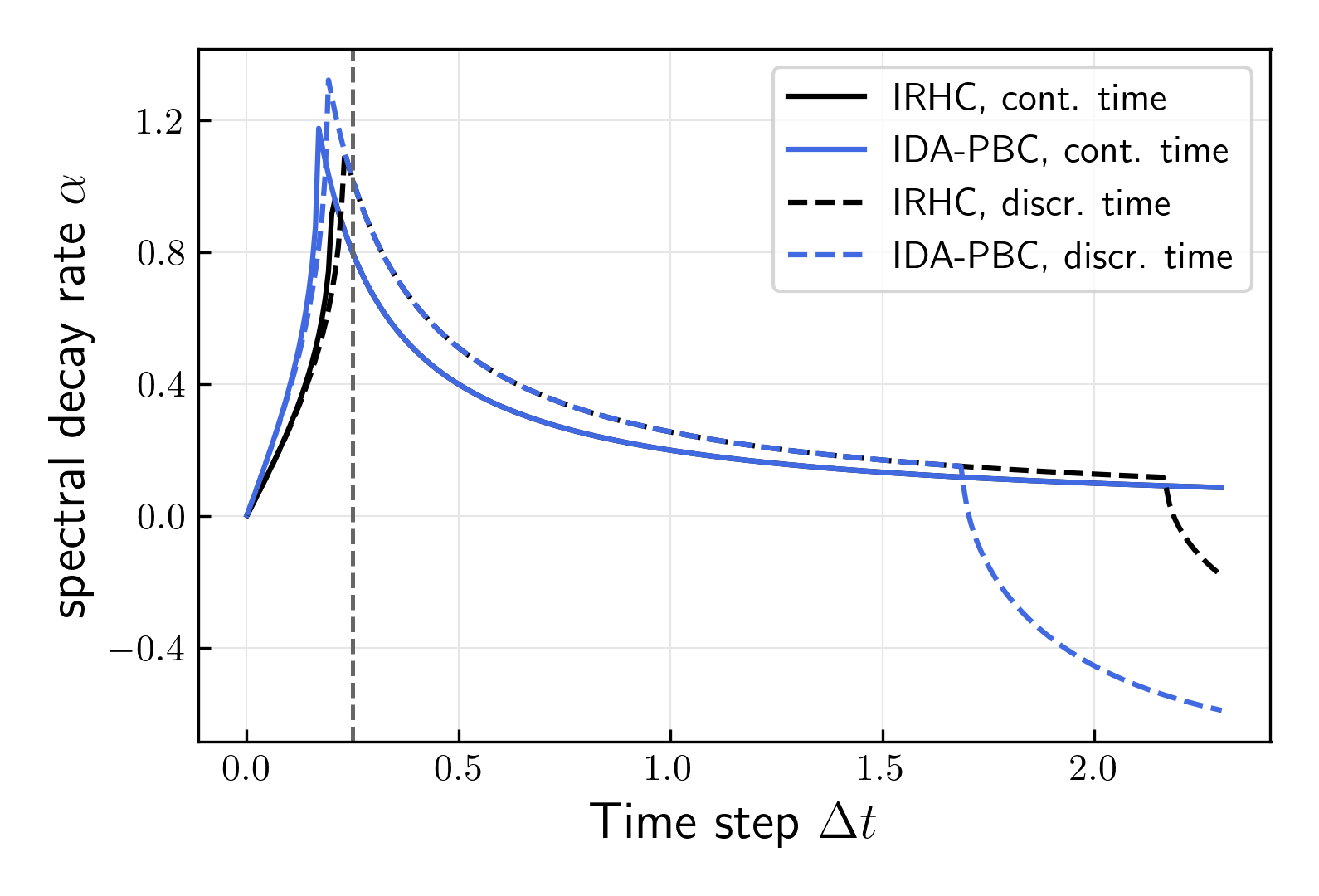}
        \caption{Discrete and continuous time spectral decay rates as a function of $\Delta t$ for the IRHC and IDA-PBC controllers. The parameters are $\B = I$, $\eta = 0.4$ and $c= 1$. The dashed vertical line marks the parameter $\Delta t = 0.25$.}
     \label{fig:agent_delta_t_decay}
    \end{figure}

    \begin{Def}
       For a square matrix $A$, we denote its spectral abscissa by $s(A)$ and its spectral radius by $r(A)$. For the discrete-time definitions, we will subsequently always assume that $r(A) > 0$ and $\Delta t > 0$. We define the continuous- and discrete-time \textit{dominant angular frequencies} of $A$ by  
       \begin{align*}
           \omega_{\text{ct}}(A) &\coloneq \max_{\substack{\lambda \in \sigma(A) \\ \mathrm{Re}(\lambda) = s (A) }} \abs{ \mathrm{Im} (\lambda) } , \\  
           \omega_{\text{dt}}(A, \Delta t) &= \frac{1}{\Delta t} \max_{\substack{\lambda \in \sigma(A) \\ \abs{\lambda} = r (A) }} 
           \abs{\arg(\lambda) },
       \end{align*}
       where $\arg(z) \in (-\pi,\pi]$ denotes the principal argument of a complex number $z$. Therefore, $\omega_{\text{ct}}$ is the largest angular frequency associated to an eigenvalue of $A$ with maximal real part. Similarly, $\omega_{\text{dt}}$ is the largest angular frequency corresponding to an eigenvalue with maximal modulus and a time step $\Delta t > 0$. Analogously, we call $\alpha_{\text{ct}} (A) \coloneq - s(A)$ the continuous-time spectral decay rate. Its discrete-time spectral decay rate is given by   
           $\alpha_{\text{dt}}(A, \Delta t) = - \ln (r(A) ) / {\Delta t }$  .
   \end{Def}
   Due to Lemma \ref{lemma:COMspectrum}, the effect of the $c_2$-term and of the time discretization on the spectral decay rate and the dominant angular frequency can be characterized explicitly. We start by introducing the distinction between over- and underdamped modes.
   \begin{Def}
       Let $\lambda \in \sigma(\B)$. We call the corresponding continuous-time mode \textit{overdamped} if $\mu_{\mathrm{ct},+}(\lambda)$ and $\mu_{\mathrm{ct},-}(\lambda)$ are real and distinct. If they form a non-real complex conjugate pair, we call the mode \textit{underdamped}. By replacing $\mu_{\text{ct}, \pm}(\lambda)$ by $\mu_{\text{dt}, \pm}(\lambda)$, we analogously define discrete-time over- and underdamped modes.
   \end{Def}
    The following lemma compares the spectral decay rates of the discrete- and continuous-time systems.
   \begin{Lemma}\label{lemma:ratecomp}
       Let $\lambda \in \sigma (\B) $ and $\Delta t > 0, 0 < \eta < 1$ and $c > 0$. Consider either the IRHC or the IDA-PBC controller. 
       \thmpart\label{lemma:ratecompI} If the discrete-time mode is underdamped, then the continuous-time mode is also underdamped. Similarly, if the continuous-time mode is overdamped, then the discrete-time mode is overdamped.
       \thmpart\label{lemma:ratecompII} If the discrete-time mode is underdamped, then 
       \begin{align*}
            \alpha_{\mathrm{dt}}(A_{\mathrm{dt}}(\lambda),\Delta t)
            > 
            \alpha_{\mathrm{ct}}(A_{\mathrm{ct}}(\lambda)).
       \end{align*}
       \thmpart\label{lemma:ratecompIII} If the continuous-time mode is overdamped, then 
       \begin{align*}
            \alpha_{\mathrm{dt}}(A_{\mathrm{dt}}(\lambda),\Delta t)
            < 
            \alpha_{\mathrm{ct}}(A_{\mathrm{ct}}(\lambda)).
       \end{align*}
   \end{Lemma}
   In summary, Lemma \ref{lemma:ratecomp} shows that the time discretization increases the spectral decay rate for discrete-time underdamped modes. For continuous-time overdamped modes, the situation is reversed. The proof of Lemma \ref{lemma:ratecomp} can be found in the Appendix. 
   \begin{figure}[ht]
        \centering
        \begin{subfigure}{\linewidth}
            \centering
            \includegraphics[width=\linewidth]{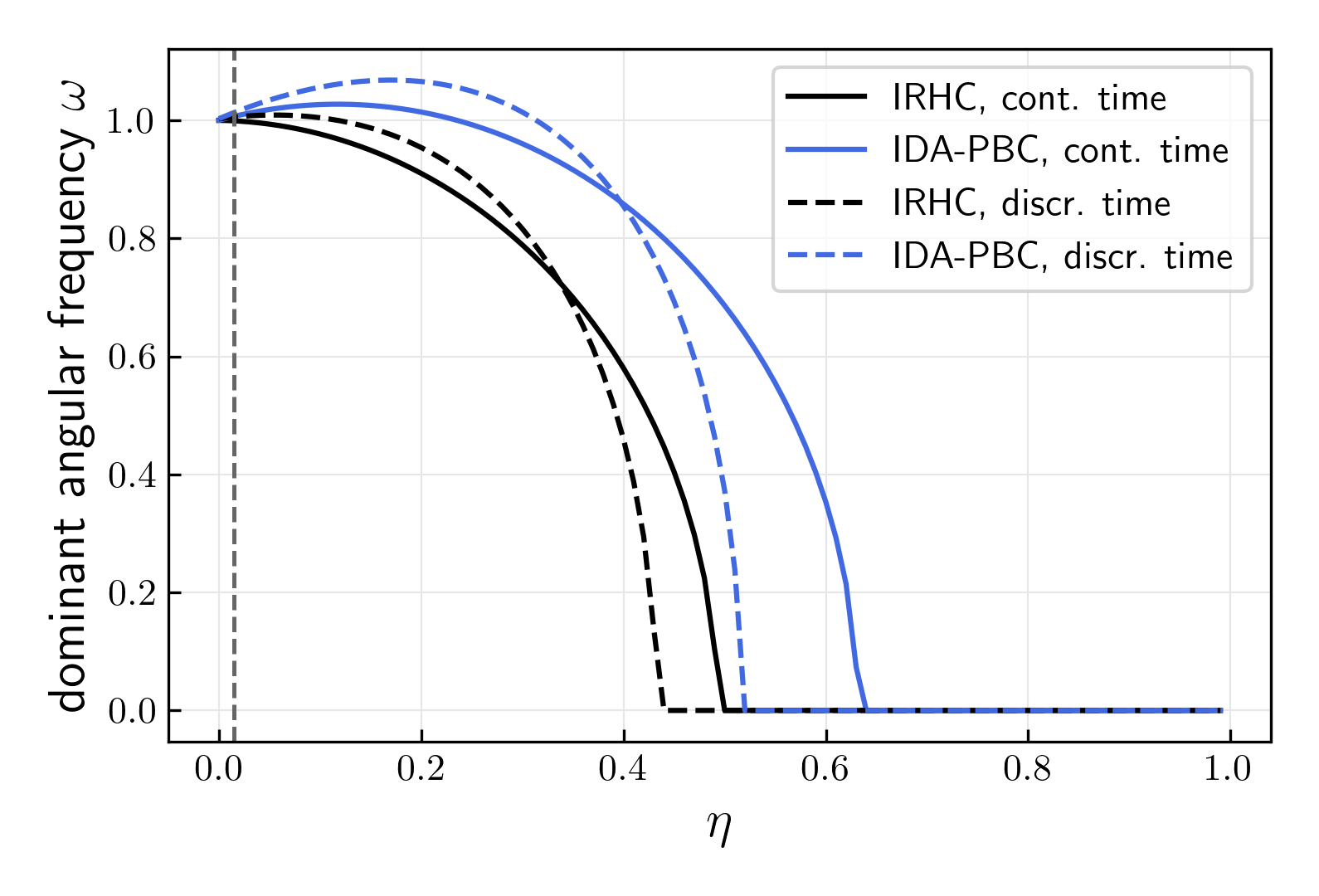}
        \caption{Discrete and continuous time dominant angular frequencies as a function of $\eta$ for the IRHC and IDA-PBC controllers.  }
        \end{subfigure}

        \vspace{0.4em}

        \begin{subfigure}{\linewidth}
            \centering
            \includegraphics[width=\linewidth]{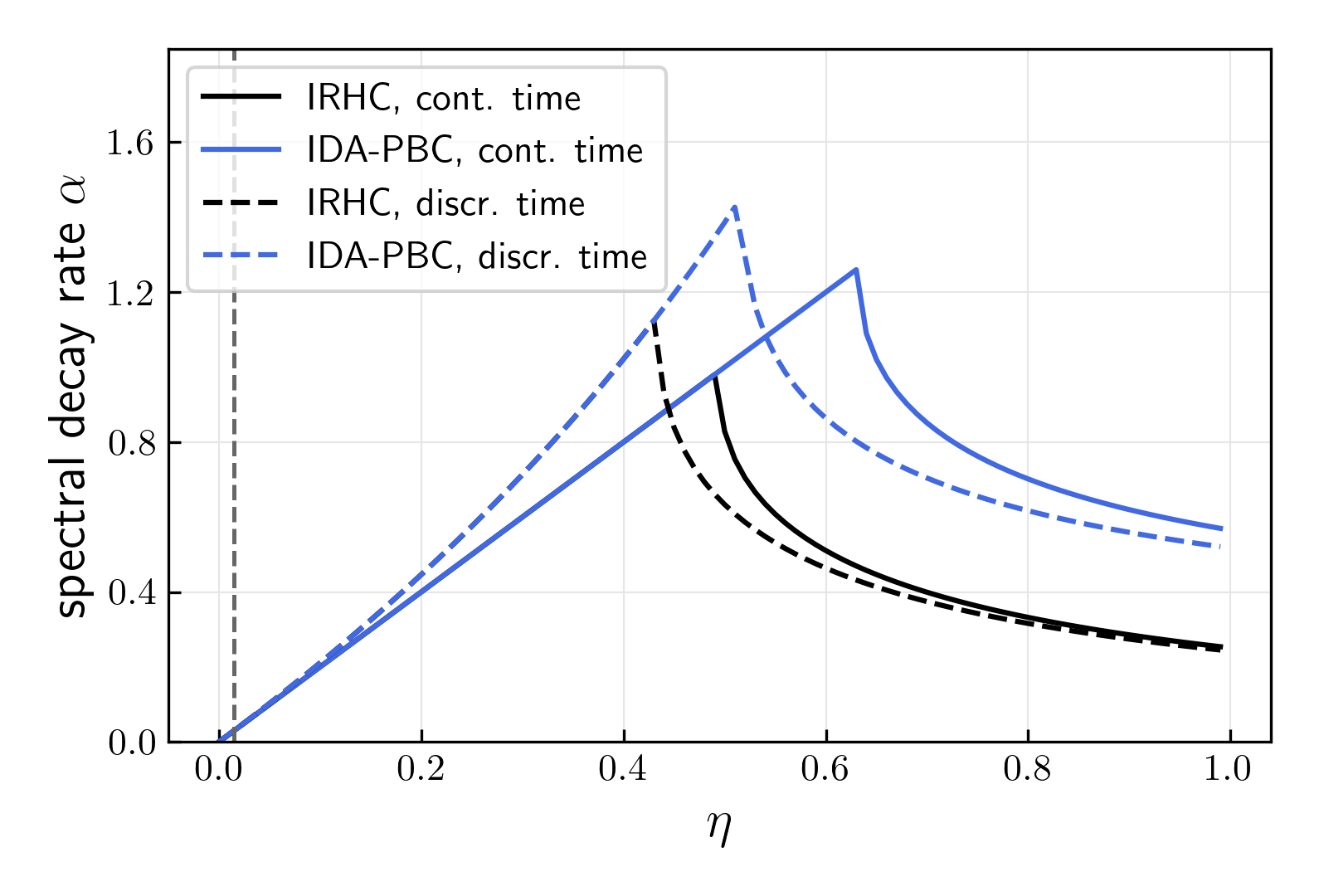}
            \caption{Discrete and continuous time spectral decay rates as a function of $\eta$ for the IRHC and IDA-PBC controllers.}
        \end{subfigure}

        \caption{The parameters are $\B = I$, $\Delta t = 0.25$ and $c= 1$. The dashed vertical lines mark the parameter $\eta = 1/{66}$ as in Theorem \ref{thm:baselinediscrtimes}.}
     \label{fig:agent_eta_decay}
    \end{figure}
   Next, we consider the effect of the parameter $c_2$ to compare the IRHC and the IDA-PBC controller. For an underdamped continuous-time mode, \eqref{eq:COMlambdaCT} implies that $ \alpha_{\text{ct}} (A_{\text{ct}  } (\lambda) ) =  c_1 \lambda^2 / {2}$, which is independent of $c_2$. Consequently, the IRHC and IDA-PBC controllers achieve the same spectral decay rate as long as the corresponding modes are underdamped for both controllers. 
   However, their angular frequencies differ. Increasing the value of $c_2$ decreases the modulus of the argument of the square root in \eqref{eq:COMlambdaCT} and therefore reduces the oscillation frequency. For an overdamped continuous-time mode, the dominant eigenvalue is $\mu_{\text{ct},+}$. This eigenvalue increases with $c_2$ and therefore the spectral decay rate decreases. Hence, overdamped modes converge faster for the IDA-PBC controller than for the IRHC controller. 
   In order to illustrate the effect of the $c_2$-term and the time discretization, we provide two figures. Fix $\B = I$ and $c = 1$. Figure \ref{fig:agent_eta_decay} shows the decay rates and oscillation frequencies for $\Delta t = 0.25$ as a function of $\eta$. The transition from under- to overdamped modes occurs at the maxima of the spectral decay rates. At those points, the corresponding dominant angular frequency vanishes. Once the continuous-time modes are overdamped, the continuous-time spectral decay rates are larger than their discrete-time counterparts. \\
   To close this section, we briefly comment on the influence of $\Delta t$. Unlike $\eta$, the parameter $\Delta t$ changes both the controller coefficients and the time step of the discrete-time closed loop. Sufficiently large values of $\Delta t$ may destabilize the discrete-time system, that is $\alpha_{\text{dt}}(A_{\text{dt} } (\lambda) , \Delta t  ) < 0  $, even though the corresponding continuous-time closed loop remains stable. The influence of $\Delta t$ on the spectral decay rate for the fixed value $\eta = 0.4$ is illustrated in Figure \ref{fig:agent_delta_t_decay}.
   Again, the transition between over- and underdamped modes occurs at the maxima of the spectral decay rate, where the overdamped regime corresponds to smaller values of $\Delta t$. For sufficiently large $\Delta t$, the discrete-time spectral decay rates in Figure \ref{fig:agent_delta_t_decay} become negative, indicating a loss of stability, whereas the corresponding continuous-time systems remain stable.

   \section{Numerical experiments}\label{sec:numerics}

   In the present section, we numerically implement the IRHC and IDA-PBC controllers for the reduced center of mass system \eqref{eq:openloopCOM} and the full nonlinear system \eqref{eq:xvparticles}. The code and the simulation data can be accessed at \cite{Zenodofiles}. 
   For the continuous-time simulations, we use the SciPy DOP853 integrator with relative and absolute tolerance values of $ \text{rtol}=  10^{-10}$ and $ \text{atol} = 10^{-11}$. For the discrete-time simulations, we use the semi-explicit Euler scheme \eqref{eq:semexupd}.
   
   \subsection{Rate and frequency verification for center of mass stabilization}
\begin{figure}[h]
        \centering
            \includegraphics[width=\linewidth]{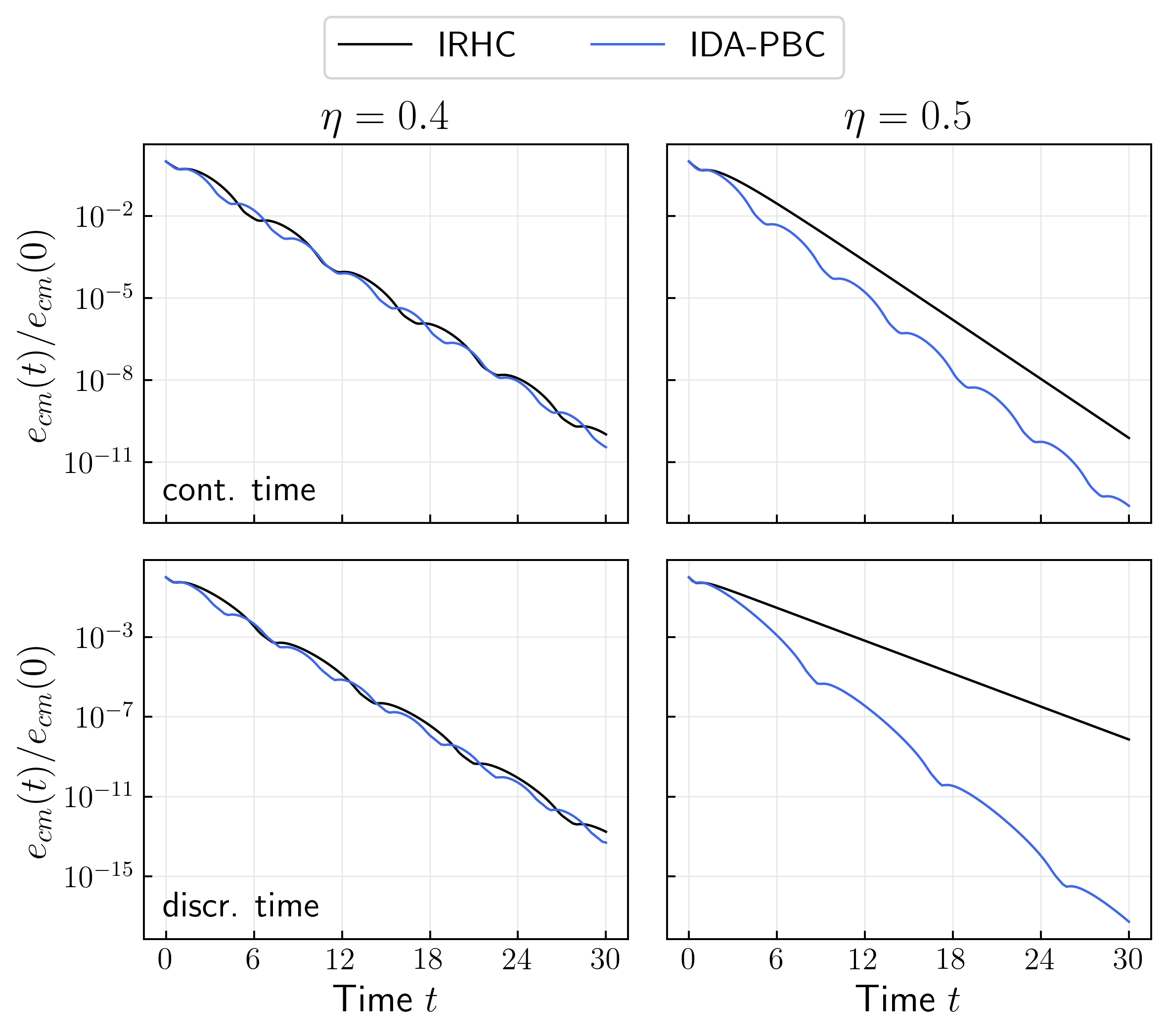}
        \caption{Simulation of the normalized center of mass error in discrete and continuous time.
        The parameters are $\B = I, c= 1$ and $\Delta t = 0.25$. The first row displays the continuous-time simulations, whereas the second row contains the discrete-time results.}
        \label{fig:agent_num_fit_decay_rates}
    \end{figure}

    We start by determining the spectral decay rates and dominant angular frequencies for the reduced center of mass system and comparing them to the predictions in Lemma \ref{lemma:COMspectrum}. For this purpose, we simulate the center of mass system in two dimensions using the parameters as in Figure \ref{fig:agent_eta_decay}, that is $\B = I$, $c= 1$ and $\Delta t = 0.25, \meanXeq = 0$, for two different values of $\eta$, namely $\eta = 0.4$ and $0.5$. Those choices are based on the predictions in Figure \ref{fig:agent_eta_decay}. For $\eta = 0.4$, both controllers are in the underdamped regime, whereas for $\eta = 0.5$ the IRHC mode is overdamped in discrete and continuous time whereas the IDA-PBC controller remains underdamped. To measure the decay of the mean velocity and the target distance of the center of mass, we use the normalized error $e_{\text{cm}}(\bar{x}(t), \bar{v}(t)) / { e_{\text{cm}}(\bar{x}(0), \bar{v}(0)) } $, where 
    \begin{align*}
        e_{\text{cm}}( \bar{\charX}, \bar{\charV})  \coloneq \abs{ \bar{\charX} - \meanXeq } + \abs{\bar{\charV} }, \ \bar{\charX}, \bar{\charV} \in \R^d.  
    \end{align*}
    Figure \ref{fig:agent_num_fit_decay_rates} shows the decay of the normalized error on a semi-logarithmic scale. 
    In Figure \ref{fig:agent_num_fit_decay_rates}, we can see that all errors exhibit an approximately linear decay in the logarithmic plot, which is consistent with the expected exponential decay. For $\eta = 0.4$, the IDA-PBC and the IRHC show similar decay rates and both errors show visible oscillations, which is consistent with the predicted underdamped regime. For $\eta = 0.5$, the IDA-PBC controller appears to have a higher spectral decay rate but its error still oscillates. No oscillations are visible for the IRHC controller, which is consistent with the expected overdamped decay. 
    In order to quantitatively compare the simulation results with the predictions of Lemma \ref{lemma:COMspectrum}, we jointly fit the center of mass position and velocity in each component. Based on the regime predicted by Lemma \ref{lemma:COMspectrum}, we choose an over- or an underdamped fitting ansatz. The fit does not use the corresponding theoretical decay rates and frequencies. For underdamped modes, we fit an exponentially damped oscillation with common decay rate $\alpha$ and angular frequency $\omega$. For overdamped modes, we fit a linear combination of two exponentials with separate decay rates and determine the spectral decay rate as the smaller of the two fitted decay rates. In continuous time, we additionally constrain the velocity fit to be the time derivative of the position fit, whereas in discrete time, the fit uses separate amplitudes for position and velocity and enforces common decay rates and a common frequency. Since $ \B= I$ gives the same modal quantities in both spatial directions, we report the median of the fitted values of the individual components. Table \ref{tab:isotropic_decay} compares the fitted values to the values predicted by Lemma \ref{lemma:COMspectrum}. The fitted spectral decay rates and, for the underdamped cases, the dominant angular frequencies agree with the predictions.

\begin{table}[h]
    \centering
    \begin{tabular}{lcccc}
        & \multicolumn{2}{c}{$\alpha$}
        & \multicolumn{2}{c}{$\omega$} \\
        \cmidrule(lr){2-3}
        \cmidrule(lr){4-5}
        Method
        & Theory & Fit
        & Theory & Fit \\
        \midrule

        \multicolumn{5}{l}{\textit{$\eta=0.4$, continuous-time}} \\
        IRHC
        & $0.8000$ & $0.8000$
        & $0.5801$ & $0.5801$ \\
        IDA-PBC
        & $0.8000$ & $0.8000$
        & $0.8582$ & $0.8582$ \\

        \addlinespace
        \multicolumn{5}{l}{\textit{$\eta=0.4$, discrete-time}} \\
        IRHC
        & $1.0217$ & $1.0217$
        & $0.4598$ & $0.4598$ \\
        IDA-PBC
        & $1.0217$ & $1.0217$
        & $0.8546$ & $0.8546$ \\

        \midrule

        \multicolumn{5}{l}{\textit{$\eta=0.5$, continuous-time}} \\
        IRHC
        & $0.8285$ & $0.8285$
        & $0.0000$ & $-$ \\
        IDA-PBC
        & $1.0000$ & $1.0000$
        & $0.6860$ & $0.6860$ \\

        \addlinespace
        \multicolumn{5}{l}{\textit{$\eta=0.5$, discrete-time}} \\
        IRHC
        & $0.6334$ & $0.6334$
        & $0.0000$ & $-$ \\
        IDA-PBC
        & $1.3863$ & $1.3863$
        & $0.3724$ & $0.3724$ \\

        \bottomrule
    \end{tabular}
    \caption{Comparison of theoretically predicted and numerically fitted spectral
    decay rates $\alpha$ and angular frequencies $\omega$.}
    \label{tab:isotropic_decay}
\end{table}

        \subsection{Full set-point stabilization}

        Next, we apply the IRHC and IDA-PBC controllers to the full nonlinear system \eqref{eq:xvparticles}. In all examples, we use the constant alignment function $\psi(x) = 1$ and the binary interactions are given by the (regularized) Morse potential 
        \begin{align}
            \mathcal{V}(x) = R \mathrm{e}^{-\abs{x}^2 / {r} } - A \mathrm{e}^{- \abs{x}^2 / a}, \label{eq:Morse}
        \end{align}
        with parameters $R = 1, r = 0.5 , A = 1.5 , a = 1$. We consider $N = 42$ particles moving in $d = 2$ dimensions. The Morse potential $\mathcal{V}$ attains its global minimum at the distance 
        \begin{align*}
            \abs{x} = d_0 = \Big( \frac{r a}{a-r}  \ln \Big( \frac{a R}{r A} \Big) \Big)^{\frac{1}{2}} \approx 0.5364.
        \end{align*}
        To obtain a target configuration $\Xeq$ containing interparticle distances on both sides of $d_0$, we sample the target position for each particle uniformly from $[-6d_0, 6d_0]^2$. The sampling is done once and the resulting configuration is then fixed for all simulations. 
        \begin{figure*}[h]
        \centering
        \includegraphics[width=\textwidth]{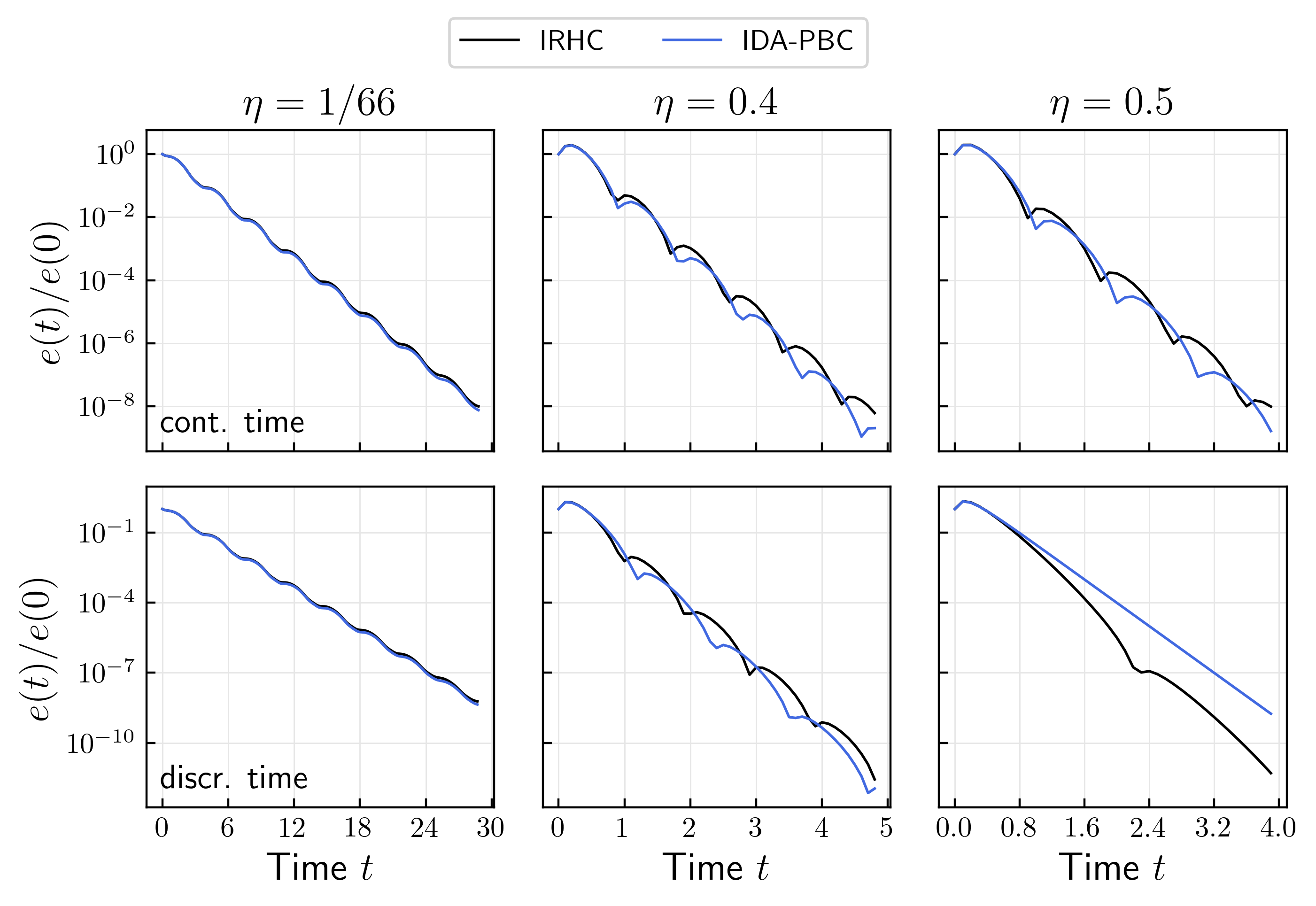}
        \caption{Simulation of the normalized centered error $e(t) / e(0) = e(\charX(t), \charV (t)) / e(\charX(0), \charV (0))$ in discrete and continuous time for the parameters $\cB = I, c= 100$ and $\Delta t = 0.05$. The first row displays the continuous-time simulations, whereas the second row contains the discrete-time results.}
        \label{fig:nonlinear_decay}
        \end{figure*}
        The aim of this subsection is to study the controller behavior for different values of $\eta$ and to relate the chosen parameters to the sufficient convergence results in Theorem \ref{thm:convbaseline}.  
        For the Morse potential \eqref{eq:Morse}, an elementary computation of the Hessian shows that $\nabla \mathcal{V} \colon \R^d \to \R^d$ is Lipschitz continuous with Lipschitz constant 
        $L = 2 ( R/{r} + A/{a} )$.
    For $\cB = I, c = 100, \Delta t = 0.05$ and the parameters $A,a,R,r$ considered here, we obtain $\sigma = 1$ and $K = 28$. Moreover, the conditions \eqref{eq:ctdtc} in Theorem \ref{thm:convbaseline} are satisfied, whereas \eqref{eq:cteta} reduces to $\eta > 0.35$. Thus the previously considered values $\eta = 0.4$ and $\eta = 0.5$ lie within the parameter range of Theorem \ref{thm:baselineconttimes}. As an additional reference, we consider $\eta = 1/{66}$, which is the upper bound appearing in the discrete-time convergence result of Theorem \ref{thm:baselinediscrtimes}.
    \begin{table}[!h]
    \centering
    \small
    \setlength{\tabcolsep}{3pt}
    \begin{tabular}{llccc}
        Method &
        & $\eta=1/66$ & $\eta=0.4$ & $\eta=0.5$ \\
        \midrule
        \multicolumn{5}{l}{\textit{Continuous time}} \\
        \addlinespace
        IRHC    & $\alpha_{\text{ct}}  $
                & 0.6297 & 4.2964 & 5.2500 \\
                & $\omega_{\text{ct}}  $
                & 0.8584 & 3.6831 & 3.5251 \\
        IDA-PBC & $\alpha_{\text{ct}}  $
                & 0.6358 & 4.4976 & 5.4903 \\
                & $\omega_{\text{ct}}  $
                & 0.8580 & 3.4159 & 3.1485 \\
        \midrule
        \multicolumn{5}{l}{\textit{Discrete time}} \\
        \addlinespace
        IRHC    & $\alpha_{\text{dt}}  $
                & 0.6522 & 5.6213 & 7.4444 \\
                & $\omega_{\text{dt}}  $
                & 0.8746 & 3.2535 & 1.4254 \\
        IDA-PBC & $\alpha_{\text{dt}}  $
                & 0.6592 & 5.9458 & 5.7239 \\
                & $\omega_{\text{dt}}  $
                & 0.9165 & 2.6364 & 0.0000 \\
        \bottomrule
    \end{tabular}
    \caption{DMD spectral decay rates and dominant angular frequencies for the centered nonlinear system.}
    \label{tab:nonlinear_decay_rates}
\end{table}
     However, we emphasize that the condition $\eta > 0.35$ is violated for $\eta = 1/{66}$, so Theorem \ref{thm:baselineconttimes} does not apply. None of the three cases for $\eta$ is covered by Theorem \ref{thm:baselinediscrtimes}. We note that for $\cB = I$, the center of mass follows a reduced dynamics. Concretely, since the binary alignment and binary interaction forces have zero mean, we obtain for both controllers
    \begin{align*}
        \dot{ \overline{\charX} }  &= \overline{\charV}, \\
        \dot{ \overline{\charV} } &= - c_1 \overline{\charV} - c_3
        (\overline{\charX} - \meanXeq),
    \end{align*}
    for which the spectrum of the system matrix can be computed similarly as in Lemma \ref{lemma:COMspectrum}. In particular, we can explicitly predict the continuous-time spectral decay rate of the center of mass and this rate is independent of $c_2$ and therefore identical for both controllers. 
    For all parameters considered here, the corresponding mode is underdamped and thus the continuous-time center of mass spectral decay rate is given by $c_1 /{2}$.
    Therefore, we isolate the internal particle dynamics from the explicitly known center of mass mode by subtracting the mean position and mean velocity. More precisely, we use the normalized centered error $e (\charX (t), \charV (t)) / { e (\charX (0), \charV (0)) } $, where 
    \begin{align*}
        e( \charX , \charV )  \coloneq \norm{ \charX - \Xeq - (\overline{\charX} - \meanXeq) \mathds{1} } + \norm{ \charV - \overline{\charV} \mathds{1} },
    \end{align*}
    for $ \charX, \charV \in \HS  $. 
    Figure \ref{fig:nonlinear_decay} shows the decay of the normalized centered error on a semi-logarithmic scale. 
    For all considered values of $\eta$, the normalized centered error appears to exhibit an approximately linear decay on the logarithmic scale, indicating exponential decay. 
    To quantify the decay and oscillation frequencies of the centered particle dynamics, we apply dynamic mode decomposition (DMD) to the centered simulation data and compute the system matrix using a least squares fit (see \cite[Chapter $7.2$]{BruntonDMD}). The spectral decay rate and dominant angular frequencies are then determined from the spectrum of the resulting matrix. 
    Since the center of mass contribution has been removed, these fitted rates are not predicted by Lemma \ref{lemma:COMspectrum}. The resulting values are displayed in Table \ref{tab:nonlinear_decay_rates}.
    For all considered values of $\eta$, we see that the discrete-time DMD spectral decay rates are larger than their continuous-time analogue. This is consistent with the result for the reduced system that time discretization increases the spectral decay rate of modes that remain underdamped. Moreover, the IDA-PBC controller achieves a larger spectral decay rate in all cases except for the discrete-time simulation with $\eta = 0.5$. Although IDA-PBC achieves faster decay in most simulations, there is no uniform ordering of the decay rates for the centered nonlinear dynamics.
    The IRHC controller exhibits higher dominant angular frequencies than the IDA-PBC controller for $\eta = 0.4$ and $\eta = 0.5$. Therefore, the reduction in oscillation frequency we observed for the reduced system does in general not carry over to the centered nonlinear dynamics. 
    This difference is most pronounced in discrete-time at $\eta = 0.5$. The IRHC controller achieves a higher estimated spectral decay rate but retains an underdamped dominant mode, whereas IDA-PBC exhibits a lower decay rate with an overdamped dominant mode.

      \section{Conclusion}  
      The reduced center of mass dynamics provide an interpretation of the parameter $c_2$ in the IRHC feedback controller. It modifies the transition between over- and underdamped modes. For discrete-time underdamped modes, it decreases the angular frequency without changing the spectral decay rate. For continuous-time overdamped modes, it reduces the spectral decay rate. The time discretization by the semi-explicit Euler scheme adds an additional mode-dependent effect. It increases the spectral decay rate for modes that remain underdamped, whereas the reverse effect occurs for continuous-time overdamped modes. We verified these predicted effects in simulations of the reduced center of mass system. For the full nonlinear system, we also observe larger DMD-estimated spectral decay rates in discrete- than in continuous-time. However, the reduction of the dominant angular frequency for the IRHC controller does in general not carry over to the nonlinear dynamics. Explicit decay rate estimates for the centered nonlinear dynamics would complement the present analysis and could provide a basis for a rigorous tuning strategy for nonlinear particle dynamics.

    \appendix 
    \section{Proof of Lemma \ref{lemma:ratecomp}}
    \begin{proof}[Proof of Lemma \ref{lemma:ratecomp}]
        Fix $\lambda \in \sigma (\B)\subseteq (0,\infty)$. We start by rewriting the expressions for $\mu_{\text{ct},\pm}(\lambda)$. Define 
        \begin{align*}
           k \coloneq \det (A_{\text{ct}} (\lambda) )
           = \lambda + c_3 \lambda^2 - c_2 \lambda^3 .
        \end{align*}
        We note that $k > 0$. Indeed, for the IDA-PBC controller, we have $c_2 = 0$, whereas for the IRHC controller, \eqref{eq:cetaPar} yields 
        \begin{align}
            c_2 \lambda^2 \leq c_2 \norm{\B^2} =  c_2 \norm{\cB \cB^*} = \eta < 1 \label{eq:c2lest}
        \end{align}
        and therefore $k > c_3 \lambda^2 \geq 0$.
         Moreover, a short computation using \eqref{eq:COMlambdaDT} shows that 
        \begin{align}
            &\mu_{\text{dt},\pm} (\lambda) 
            -1
            \nonumber \\
            &= 
            \Delta t \Big( 
            - \frac{c_1 \lambda^2 + \Delta t k }{2} 
            \pm 
            \sqrt{ 
            \Big(\frac{c_1 \lambda^2 + \Delta t k }{2} \Big)^2 
            - k
            }
            \Big) \label{eq:COMdtrewrite}.
        \end{align}
        In particular, \eqref{eq:COMdtrewrite} and $k > 0$ imply that the implications of Lemma \ref{lemma:ratecompI} hold. Indeed, the arguments of the square root in \eqref{eq:COMlambdaCT} and \eqref{eq:COMdtrewrite} differ by a positive term. \\ 
        To show Lemma \ref{lemma:ratecompII}, let us assume that the discrete-time mode is underdamped. By Lemma \ref{lemma:ratecompI}, the corresponding continuous-time mode is also underdamped. In particular, we have $\alpha_{\text{ct}} (A_{\text{ct}}(\lambda) )  = c_1 \lambda^2/{2}$ and 
        \begin{align*}
            \abs{\mu_{\text{dt}, \pm }(\lambda) } = \sqrt{ \det (A_{\text{dt}} (\lambda) ) }
            = \sqrt{1 - c_1 \lambda^2 \Delta t} .
        \end{align*}
        As in \eqref{eq:c2lest}, we note that $c_1 \Delta t \lambda^2 < 1$. Thus, it follows that 
        \begin{align*}
            \alpha_{\text{dt} } &= 
             - \frac{1}{2 \Delta t} 
            \ln ( 1 - c_1 \lambda^2 \Delta t )
            > \frac{c_1 \lambda^2}{2} = \alpha_{\text{ct}} ,
        \end{align*}
        where we have used that $- \ln(1-x) > x$ for $x \in (0,1)$. \\
        Now let us assume that the continuous-time mode (and thus also the discrete-time mode) is overdamped. From \eqref{eq:COMlambdaCT}, we see that 
            $k < (c_1 \lambda^2 / {2} )^2 $
        and therefore 
        \begin{align*}
            \mathrm{tr}(A_{\text{dt} }(\lambda)  )
            &= 2 - c_1 \Delta t \lambda^2 
            - (\Delta t)^2 k \\
            &> 2 - c_1 \Delta t \lambda^2  
            -
            \frac{(c_1 \lambda^2 )^2 (\Delta t)^2}{4} \\
            &\geq 2 - \eta - \frac{\eta^2}{4} > 0.
        \end{align*}
        Thus we conclude that $r(A_{\text{dt}}(\lambda) ) = \mu_{\text{dt}, + } (\lambda) > 0$.
        Since the continuous-time mode is overdamped, we also note that $\alpha_{\text{ct}} (A_{\text{ct} }(\lambda ) ) = -\mu_{\text{ct}, + }(\lambda)$. 
        We define the function 
        \begin{align*}
            g \colon (\sqrt{k}, \infty) \to \R, \ g(s) = - s + \sqrt{s^2 - k}.
        \end{align*}
        We collect some elementary facts about the function $g$. First, $g$ is continuously differentiable and strictly increasing. Moreover, we have $\range(g) = (-\sqrt{k},0)$ and the inverse of $g$ is given by 
        \begin{align}
            g^{-1} (y) =  - \frac{k}{2y} - \frac{y}{2} , \ y \in (-\sqrt{k}, 0) . \label{eq:invg}
        \end{align}
        We claim that the following inequality holds 
        \begin{align}
            \frac{c_1 \lambda^2 + k \Delta t}{2} > 
            g^{-1} \Big( 
            \frac{\mu_{\text{ct},+ }(\lambda) }{1 - \Delta t \mu_{\text{ct},+ }(\lambda) } 
            \Big). \label{eq:mainestimate}
        \end{align}
        In particular, we claim that the argument of $g^{-1}$ above lies in $(-\sqrt{k},0)$.
        Once \eqref{eq:mainestimate} is shown, the proof is complete because then we have 
        \begin{align*}
            \mu_{\text{dt},+}(\lambda) 
            &= 1 + \Delta t g \Big( \frac{c_1 \lambda^2 + k \Delta t}{2} \Big) \\
            &> 1 + \Delta t \frac{\mu_{\text{ct},+ }(\lambda) }{1 - \Delta t \mu_{\text{ct},+ }(\lambda) } \\
            &= \frac{1}{1- \Delta t \mu_{\text{ct},+ } (\lambda) } 
            \geq \mathrm{e}^{\Delta t \mu_{\text{ct}, + }(\lambda) } .
        \end{align*}
        Taking logarithms and dividing by $\Delta t$ then yields 
        \begin{align*}
            \alpha_{\text{dt}} (A_{\text{dt} } (\lambda), \Delta t ) 
            &= - \frac{\ln( \mu_{\text{dt},+}(\lambda) )}{\Delta t}  < -\mu_{\text{ct},+}(\lambda) \\
            &= \alpha_{\text{ct}}(A_{\text{ct}} (\lambda) ) .
        \end{align*}
        So to complete the proof, it remains to show \eqref{eq:mainestimate}. For this purpose, we first note that \eqref{eq:mainestimate} is meaningful. Indeed, since $\mu_{\text{ct},+}(\lambda) < 0$ and since the continuous-time mode is overdamped, we have $\abs{\mu_{\text{ct}, - } (\lambda) } > \abs{\mu_{\text{ct}, + }(\lambda) }$. Thus 
        \begin{align*}
            k  = \abs{  \mu_{\text{ct}, - } (\lambda)    \mu_{\text{ct}, + }(\lambda) } 
            > \abs{ \mu_{\text{ct}, + }(\lambda) }^2
        \end{align*}
        and therefore 
        \begin{align*}
            0 > \frac{\mu_{\text{ct},+ }(\lambda) }{1 - \Delta t \mu_{\text{ct},+ }(\lambda) }
            > \mu_{\text{ct},+ }(\lambda) > - \sqrt{k}.
        \end{align*}
        By using \eqref{eq:COMlambdaCT}, we see that $\mu_{\text{ct},+}(\lambda)$ satisfies the quadratic equation 
        \begin{align}
            0 = \mu_{\text{ct}, +}(\lambda)^2 + c_1 \lambda^2 \mu_{\text{ct}, +}(\lambda) + k. \label{eq:mu+poleq}
        \end{align}
        By using \eqref{eq:mu+poleq}, a short computation shows that 
        \begin{align}
            &- \frac{k}{\mu_{\text{ct}, + }(\lambda) } (1 - \Delta t \mu_{\text{ct}, + }(\lambda) ) \nonumber \\
            &= \mu_{\text{ct}, +}(\lambda) + k \Delta t + c_1 \lambda^2. \label{eq:mainestpart}
        \end{align}
        Finally, combining \eqref{eq:invg} and \eqref{eq:mainestpart} yields
        \begin{align*}
            &g^{-1} \Big( 
            \frac{\mu_{\text{ct},+ }(\lambda) }{1 - \Delta t \mu_{\text{ct},+ }(\lambda) } 
            \Big) \\
            &= \frac{\mu_{\text{ct}, +}(\lambda) + k \Delta t + c_1 \lambda^2}{2} 
            -  \frac{\mu_{\text{ct},+ }(\lambda) }{2(1 - \Delta t \mu_{\text{ct},+ }(\lambda)) } 
              \\
            &= -\frac{\Delta t \mu_{\text{ct}, + }(\lambda)^2 }{2 (1 - \Delta t \mu_{\text{ct},+ }(\lambda) )}
            + \frac{c_1 \lambda^2 + k \Delta t }{2} \\ 
            &< \frac{c_1 \lambda^2 + k \Delta t}{2},
        \end{align*}
        which completes the proof. 
    \end{proof}

    \section*{Acknowledgments}
    This work was funded by the Deutsche Forschungsgemeinschaft (DFG, German Research Foundation) – Project-ID 531152215 – CRC 1701. We acknowledge the use of ChatGPT for stylistic editing, manuscript review and assistance with refining and debugging parts of the numerical simulation code used for the experiments in Section \ref{sec:numerics}. The authors assume responsibility for all content.

\bibliographystyle{unsrtnat}
\bibliography{references}

\end{document}